\documentclass[a4paper,11pt]{amsart}
\usepackage{geometry} 
\usepackage{amsmath,amssymb,mathrsfs}
\usepackage{xcolor}
\usepackage{indentfirst}
\usepackage[colorlinks,
            linkcolor=black,
            anchorcolor=black,
            citecolor=black
            ]{hyperref}
\DeclareMathOperator{\diver}{\mathrm{div}}
\def\T{{\mathbb T}^{3}}
\def\R{{\mathbb R}^{3}}

\def\N{\mathbb N}
\def\K{\mathbb K}
\def\d{\partial_t}
\def\eps{\epsilon}
\def\O{\Omega}
\newtheorem{mydef}{Definition}
\newtheorem{thm}{Theorem}
\newtheorem{lem}{Lemma}
\newtheorem*{pf}{Proof}
\newtheorem{rk}{Remark}

\numberwithin{thm}{section}
\numberwithin{mydef}{section}
\numberwithin{lem}{section}
\numberwithin{rk}{section}
\numberwithin{equation}{section}
\title[Quantum MHD equations]{Global existence of weak solutions for quantum MHD equations}
\author{Hao Li \quad Yachun Li}
\date{}
\address[Hao.Li]{School of Mathematical Sciences, Shanghai Jiao Tong University,Shanghai 200240, P.R.China;}\email{\tt haoli\_sjtu@163.com}
\begin{document}
\maketitle
\begin{abstract}
In this paper, we consider the quantum MHD equations with both the viscosity coefficient and the magnetic diffusion coefficient are depend on the density. we prove the global existence of weak solutions and perform the lower planck limit in a 3-dimensional torus for large initial data.  The global existence is shown by using Faedo-Galerkin method and weak compactness techniques for the adiabatic exponent $\gamma>1$. 
\end{abstract}
\section{introduction}
Quantum fluid models are used to describe, for instance, superfluids \cite{Loffredo1993On}, quantum semiconductors\cite{PhysRevB.48.7944}, thermistor theory and weakly interacting Bose gases \cite{Grant1973Pressure}, and quantum trajectories of Bohmian mechanics \cite{Wyatt2006Quantum}. In this paper we mainly study the quantum MHD equations in $\O=\T$ ($\T$ is the 3-dimensional torus in $\R$) which read as follows:
\begin{equation}\label{e1.1}
\left\{
\begin{aligned}
&\d \rho+\diver(\rho u)=0, x\in \O, t>0,\\
&\d(\rho u)+\diver(\rho u\otimes u)+\nabla(P(\rho)+P_c(\rho))-2\diver(\rho D(u))\\
&\quad -2\kappa^2\rho\nabla\left(\frac{\Delta\sqrt\rho}{\sqrt\rho}\right)-(\nabla\times B)\times B=0,\\
&\d B-\nabla\times(u\times B)+\nabla(\nu_b(\rho)\nabla\times B)=0,\\
&\diver B=0.
\end{aligned}
\right.
\end{equation}
with the initial data 
\begin{equation}
\rho(0,x)=\rho_0(x), (\rho u)(0,x)=m_0,B(0,x)=B_0(x), \diver B_0=0.
\end{equation}
where the functions $\rho,u$ and $B$ represent the mass density, the velocity field and the magnetic field respectively. The function $P(\rho)=\rho^{\gamma}$ with $\gamma>1$ is the pressure, $P_c(\rho)$ is a singular continuous functions and called the cold pressure. $D(u)=\frac{\nabla u+(\nabla u)^{\top}}{2}$ is the stress tensor,  $\nu_b(\rho)$ is the magnetic diffusion viscosity coefficient, $\kappa>0$ is the quantum planck constant. The expression $2\kappa \rho\nabla\left({\Delta\sqrt \rho}/{\sqrt \rho}\right)$ can be interpreted as a quantum Bohm potential, and has the following identity:
\begin{equation} 2\kappa^2 \rho\nabla\left(\frac{\Delta\sqrt \rho}{\sqrt \rho}\right)=\kappa^2\diver(\rho\nabla^2\log \rho)=\kappa^2\nabla\Delta \rho-4\kappa^2\diver(\nabla\sqrt \rho\otimes\nabla\sqrt \rho).\end{equation}

Recently, quantum fluid models have received a great deal of attention of mathematicians.  J\"{u}ngel\cite{JDissipative} derived the dissipative quantum fluid models from Wigner-Boltzmann  equation by a moment method, and the quantum ideal magnetohydrodynamic model was derived by Hass\cite{Haas2011Quantum}. The existence of global weak solutions have been studied by many authors. For the compressible Navier-Stokes equations with constant viscosity coefficients, the pioneer work is P.L.Lions \cite{Lions1998Mathematical} who proved the global existence of weak solutions for the barotropic compressible Navier-Stokes systems with $\gamma>3n/(n+2)$. Later, Feireisl \cite{Feireisl2001} extend this result to the case $\gamma>n/2$. When the viscosity coefficients $\mu,\lambda$ are density-dependent, the systems become much more difficult due to the velocity cannot be defined in the vacuum region. Under the assumption on the viscosity coefficients i.e $\lambda(\rho)=2(\rho\mu^{'}(\rho)-\mu(\rho))$. Bresch-Desjardins \cite{BRESCH2006362,Bresch2003Existence,Bresch2004Some,Didier2003On} made great progress, they discover a new mathematical entropy inequality which is not only  applied to the vacuum case but also  applied to get the existence of global weak solutions. Mellet-Vasseur \cite{doi:10.1080/03605300600857079} study the stability of the baratropic compressible Navier-Stokes equations. Later, Vasseur-Yu \cite{Vasseur2016} and Li-Xin \cite{Li2015Global} independently prove the global weak solution for $3D$ degenerate compressible Navier-Stokes equation, where they constructed separately appropriate approximation by  different approaches.

 Now we recall some results on the compressible quantum Navier-Stokes equations. J\"{u}ngel \cite{doi:10.1137/090776068} proved the  global existence of weak solutions by choosing $\rho\varphi$ as test function. However, this particular choice of test function does not contain the region  $\{\rho(x,t)=0\}$ in the weak formulation.  Gisclon-Violet \cite{GISCLON2015106} also showed the global existence results but in the classical definition of weak solutions by adding the cold pressure term. They also pointed that the cold pressure term can be replaced by the drag force term. Vasseur-Yu \cite{doi:10.1137/15M1013730} also considered the compressible quantum Navier-Stokes equations with damping term which is helpful to get the result \cite{Vasseur2016}. 

So far there are few results on the global existence of the quantum MHD equations. Yang-Ju \cite{doi:10.1063/1.4891492} proved the existence of global weak solutions for a special parabolic systems on the density $\rho$ and on the momentum $\rho u$ by doing a transformation for the velocity. Very recently, Guo-Xie \cite{Boling2017GLOBAL} established the global existence of weak solutions for the $2D$ general quantum MHD equations for which both viscosity coefficients and the dispersion term are general function. In the present paper, we prove the global existence of weak solutions for $3D$ quantum MHD models under special assumption on the viscosity coefficients and the dispersion term. It should be noted that we also require the magnetic diffusion coefficient satisfies the specific condition which is important to get the BD entropy. In addition, we perform the lower planck limit.

In the present paper we make some assumptions that have a physical background which is similar to \cite{Boling2017GLOBAL}

\textbf{Assumption 1.} $\mu,\lambda$ are respectively the shear and bulk viscosity coefficients and satisfy the BD entropy relationship i.e $\lambda(\rho)=2(\rho\mu^{'}(\rho)-\mu(\rho))$. In this paper,we deal with  $\mu(\rho)=\rho, \lambda(\rho)=0$.

\textbf{Assumption 2.} The cold pressure $P_c(\rho)$ is a suitable increasing function satisfying 
\begin{equation} \lim_{\rho\to 0} P_c(\rho)=\infty . \end{equation}
More precisely, we assume 
\begin{equation}
P_c^{'}(\rho)=\left\{
\begin{aligned}
&c_1\rho^{-\gamma^{-}-1},\rho\leq 1,\\
&c_2\rho^{\gamma-1}, \rho>1,
\end{aligned}
\right.
\end{equation}
where $\gamma^{-},\gamma\ge 1, c_1, c_2>0$.

\textbf{Assumption 3.} The magnetic diffusion viscosity coefficient $\nu_b(\rho)$ is a continuous functions of the density, bounded from above and take large values for both small and large densities. Furthermore, we assume that there exists $M>0$, positive constants $d_0,d_1,d_2,d_3$ large enough, and $2\leq a<a^{'}<3$ such that
\begin{equation}
\left\{
\begin{aligned}
& \frac{d_0}{s^a}\leq \nu_b(s)\leq \frac{d_1}{s^{a^{'}}}, s<M,\\
&d_2\leq \nu_b(s)\leq d_2 s^b, s\ge M.
\end{aligned}
\right.
\end{equation} 
\textbf{Assumption 4.}  Functions $H(\rho)$ and $H_c(\rho)$ are satisfy the following relationship:
\begin{equation}\label{e1.7}
\rho H^{'}(\rho)-H(\rho)=P(\rho), \rho H_c^{'}(\rho)-H_c(\rho)=P_c(\rho).
\end{equation}

Our paper is organized as follows: In section 2 we collect some elementary facts and some important inequalities which will be used in the proof of our result. In section 3 we state our main result. In section 4, 5 we prove the global existence of weak solution for the approximate systems by using Faedo-Galerkin method. In section 6 we devoted to deriving the B-D entropy which plays an important role to perform the limit for the parameters. In section 7 we justify the vanish lower planck limit.

\section{preliminaries}
In this section we recall some known facts and inequalities which  will be frequently used through out the paper.
The following well-known Gargliardo-Nirenberg inequality which can be found in \cite{Zeidler1990Nonlinear}
\begin{lem}
Let $\O\subset \mathbb R^n$ be a bounded open set with $\partial\O\in C^{0,1}$, $m\in N, 1\leq p,q,r\leq \infty$. Then there exists a constant $C>0$ such that for all $u\in W^{m,p}\cap L^q$
$$\left\|D^{\alpha} u\right\|_{L^r}\leq C\left\|u\right\|_{W^{m,p}}^{\theta}\left\|u\right\|_{L^q}^{1-\theta},$$
where $0\leq \left|\alpha\right|\leq m-1,$ and $\frac{1}{r}-\frac{\alpha}{d}=\theta(\frac{1}{p}-\frac{m}{d})+(1-\theta)\frac1q$. If $m-\left|\alpha\right|-\frac{d}{p}\notin N_0$, then $\theta\in[\left|\alpha\right|/m,1]$ is allowed.
\end{lem}

The following two lemmas  will be used to get the strong convergence of the solutions through out this paper.
\begin{lem}[Aubin-Lions lemma\cite{Simon198765}]
Let $X_0,X$ and $X_1$ be three Banach space with $X_0\subset X\subset X_1$. Suppose that $X_0$ is compactly embeded in $X$ and $X$ is continuously embeded in $X_1$.
\begin{enumerate}
\item Let $G$ be bounded in $L^p(0,T;X_0)$ where $1\leq p<\infty$ and $\frac{dG}{dt}$ be bounded in $L^1(0,T;X_1)$, Then $G$ is relatively compact in $L^p(0,T;X)$.
\item Let $F$ be bounded in $L^\infty(0,T;X_0)$ and $\frac{dF}{dt}$ be bounded in $L^p(0,T;X_1)$ with $p>1$ then $F$ is relatively compact in $C([0,T];X)$.
\end{enumerate}
\end{lem}

\begin{lem}
Let $\K$  be a compact subset of $\mathbb R^n(n\ge1)$. And a sequence $v^{\eps}$ satisfy
\begin{enumerate}
\item $v^{\eps}$ is uniformly bounded in $L^{1+\alpha}(\K)$ with $\alpha>0$,
\item $v^{\eps}$ converge almost everywhere to $v$,
\end{enumerate}
then $v^{\eps}$ converge strongly to $v$ in $L^1(\K)$ with $v\in L^{1+\alpha}(\K)$.
\end{lem}
\begin{lem}[\cite{doi:10.1137/090776068,doi:10.1137/15M1013730}]
For any smooth positive function $\rho(x)$, we have
\begin{equation*}
C_1\int \left|\nabla^2\sqrt\rho\right|^2 dx+C_2\int \left|\nabla\rho^{1/4}\right|^4 dx\leq \int \rho\left|\nabla\log\rho\right|^2 dx,
\end{equation*}
where $C_1,C_2$ are positive constant.
\end{lem}
\section{main results}
In this section we present two results. The first one gives the existence of weak solutions to \eqref{e1.1} without any assumption on $\gamma$ for 3-dimensional case. The second one is devoted to lower planck limit and also shows that global weak solution of \eqref{e1.1} tends to the weak solution of \eqref{e1.1} with $\kappa=0$.
Next, we will give the definition of the weak solution to \eqref{e1.1}.

\begin{mydef}\label{e2.1}
Functions $(\rho,u,B)$ are called a weak solution to \eqref{e1.1} if the following conditions are satisfied:
\begin{enumerate}
\item The continuity equation holds in the sense of distributions, i.e 
\begin{equation}
\int\rho_0\varphi(0)+\iint(\rho\varphi_t+\sqrt\rho\sqrt\rho u\nabla\varphi)dxdt=0,
\end{equation}
for any smooth test function with compactly supported $\varphi$ such that $\varphi(T,.)=0$.
\item The momentum equation satisfies 
\begin{equation}
\begin{aligned}
&\int m_0\varphi(0)+\iint\left(\sqrt\rho(\sqrt\rho  u)\varphi_t+\sqrt\rho u\otimes\sqrt\rho u\nabla\varphi+P(\rho)\diver\varphi+P_c(\rho)\diver\varphi \right)dxdt\\
&-\iint [2(\sqrt\rho u\otimes\nabla\sqrt\rho)\nabla\varphi)\nabla\varphi-2(\nabla\sqrt\rho\otimes\sqrt\rho u)\nabla\varphi-\sqrt\rho\sqrt\rho u\Delta\varphi -\sqrt\rho\sqrt\rho u\nabla\diver\varphi ]dxdt\\
&-4\kappa^2\iint(\nabla\sqrt\rho\otimes\nabla\sqrt\rho)\nabla\varphi-2\kappa^2\iint\sqrt\rho\nabla\sqrt\rho\nabla\diver\varphi-\iint(\nabla\times B)\times B\cdot\varphi dxdt=0,
\end{aligned}
\end{equation}
for any smooth test function with compactly supported $\varphi$ such that $\varphi(T,.)=0$.
\item The magnetic field $B$ satisfies 
\begin{equation}
\int B_0\varphi(0)=\iint \left(B\varphi_t+(u\times B)\cdot(\nabla\times\varphi)-\nu_b(\rho)\nabla\times B:\nabla\varphi\right) dxdt,
\end{equation}
for any smooth test function with compactly supported $\varphi$ such that $\varphi(T,.)=0$.

 Our main results on the weak solutions reads as follows:
\begin{thm}
Assume $T>0, \gamma^{-}\ge 4, \gamma>1$, and let $(\rho_0,u_0,B_0)$ satisfies $\rho_0\ge 0$ and 
\begin{equation}
\int\left(\frac{\left|m_0\right|^2}{2\rho_0}+H(\rho_0)+H_c(\rho_0)+2\kappa^2\left|\nabla\sqrt\rho_0\right|^2+\left|B_0\right|^2\right)dx\leq C.
\end{equation}
Then, there exists a global weak solution to the problem \eqref{e1.1}-\eqref{e1.7} in the sense of distribution of Definition 3.1. In particular, the weak solution $(\rho,u,B)$ satisfies the energy estimate \eqref{e3.14} and entropy inequality \eqref{e5.2}, \eqref{e5.20}.
\begin{align*}
&\rho\ge 0, \sqrt\rho\in L^\infty(0,T;H^1)\cap L^2(0,T;H^2),\\
&\rho\in L^\infty(0,T;L^\gamma), \rho^{-1}\in L^\infty(0,T;L^{\gamma-}),\rho^\gamma\in L^{5/3}(0,T;L^{5/3}),\\
&\sqrt\rho u\in L^\infty(0,T;L^2), \nabla\left(\frac{1}{\sqrt n}\right)\in L^2(0,T;L^2),\\
&B\in L^\infty(0,T;L^2)\cap L^2(0,T;H^1).
\end{align*}
\end{thm}
\begin{rk}
It should be noted that the Assumption 3 on the magnetic diffusion coefficient is required, which plays an important role to obtain the B-D entropy inequality.
\end{rk}
\begin{rk}
Compared with Yang-Ju \cite{doi:10.1063/1.4891492} work, we make an improvement in the present paper. In \cite{doi:10.1063/1.4891492}  both the continuity equation and the momentum equation are become parabolic with respect to $\rho$ and $\rho u$ by using a transformation. 
\end{rk}
\begin{thm}
Assume $T>0, \gamma^{-}\ge 4, \gamma>1$, and  the initial data $(\rho_0,u_0,B_0)$ satisfies the assumption in Theorem 2.1.  If we assume $(\rho_\kappa,u_\kappa,B_\kappa)$  are solutions of system \eqref{e1.1} .We have when $\kappa\to 0$ the limit function $(\rho,u,B)$ is the weak solution to the problem \eqref{e1.1}-\eqref{e1.7} with $\kappa=0$.
\end{thm}
\end{enumerate}
\end{mydef}
\section{faedo-galerkin approximation}
In this section we proved the existence of approximate solutions to the systems \eqref{e1.1} by the Faedo-Galerkin method. Motivated by the work of Feireisl, Novotny and Petzeoltova, we proceed similarly as in \cite{Feireisl2009Singular,Feireisl2004Dynamics}.
\subsection{Local solvability of the approximate system }

This section is dedicated to prove the local existence of the approximate system. We adopt the following strategy:
\begin{itemize}
\item For given $u\in C([0,T];X_n) $, the approximate continuity equation can be solved directly by means of the classical theory of parabolic equations $\rho=S(u)$.
\item Given $u$, the magnetic equation is also a linear parabolic equation and we can also find a solution $B=G(u)$ by the standard Galerkin method.
\item Having solved the $\rho, B$, we can treat the approximate momentum equation as a nonlinear integral equation. The solution of the approximate equation is based on the fixed point argument in the Banach space $C([0,T];X_n)$.
\end{itemize}
We introduce the finite dimensional space $X_n=span\{e_1,e_2,\dots, e_n\}, n\in \N$,  each $e_i$ is an orthonormal basis of $L^2$ which is also the orthogonal basis of $H^2$ we notice that $u\in C([0,T];X_n)$ is given by 
\begin{equation}
u(t,x)=\sum_{i=1}^n \lambda_i(t)e_i(x), (t,x)\in [0,T]\times\O.
\end{equation}
for some functions $\lambda_i(t)$, and the norm of $u\in C([0,T];X_n)$ can be define as 
\begin{equation*}
\left\|u(x,t)\right\|_{ C([0,T];X_n)}=\sup_{t\in[0,T]}\sum_{i=1}^n \left|\lambda_i(t)\right|.
\end{equation*}
And, $u$ can be bounded in $ C([0,T];C^k)$ for any $k\ge0$, thus we have
\begin{equation}
\left\|u\right\|_{C([0,T];C^k)}\leq C(k)\left\|u\right\|_{C([0,T];L^2)}.
\end{equation}
\begin{enumerate}
\item Approximate continuity equation
\begin{equation}\label{e3.3}
\d\rho+\diver(\rho u)=\eps\Delta\rho,
\end{equation}
with the initial data\begin{equation}\rho(x,0)=\rho_0(x)\ge \mu>0, \rho_0(x)\in C^\infty.\end{equation}
where $\mu>0$ is a constant.
For given $u\in C([0,T];X_n)$, there exists a classical solution to \eqref{e3.3}. By the maximal principle we know that the density satisfies the follows inequality
\begin{equation}
\inf_{x\in\O}\rho_0(x)\exp^{-\int_0^T\left\|\diver u\right\|_{L^\infty} ds}\leq \rho(x,t)\leq \sup_{x\in \O}\rho_0(x)\exp^{\int_0^T\left\|\diver u\right\|_{L^\infty} ds}.
\end{equation}
for all $(x,t)\in[0,T]\times\O$. Furthermore, we can also get that there exist a constant $\bar\rho$ such that 
\begin{equation}
0<\bar\rho\leq \rho(x,t)\leq \frac{1}{\bar\rho}, (x,t)\in [0,T]\times\O.
\end{equation}
Thus, we can introduce a operator $S$ from $ C([0,T];X_n)$ to $C([0,T];C^k)$ by $S(u)=\rho$. and the operator is Lipschitz continuous in the following sense:
\begin{equation}
\left\|S(u_1)-S(u_2)\right\|_{C([0,T];C^k)}\leq C(n,k)\left\|u_1-u_2\right\|_{C([0,T];L^2)}.
\end{equation}
Since for given $u\in C([0,T];X_n)$, the density is solved in terms of $u$. The magnetic equation become a linear parabolic-type equation. we can find a unique solution $B\in L^\infty(0,T;L^2)\cap L^2(0,T;H^1)$ by the standard Galerkin method and satisfies the following 
\begin{equation}\label{e3.8}
\left\{
\begin{aligned}
&\d B-\nabla\times(u\times B)+\nabla\times(\nu_b(\rho)\nabla\times B)=0,\\
&\diver B=0,\\
&B(0,x)=B_0(x).
\end{aligned}
\right.
\end{equation}
In fact, if we assume $B=B_1-B_2$  where $B_1,B_2$ are two solutions of equation with the same data then we know that $B$ also satisfied the \eqref{e3.8}.  Multiplying the equation $\eqref{e3.8}_1$ by $B$ and integrate over $\O$ we get
\begin{equation}\label{e3.9}
\begin{aligned}
&\frac12\frac{d}{dt}\int \left|B\right|^2 dx+\int\nu_b(\rho)\left|\nabla\times B\right|^2 dx\\
&=\int (u\times B)\cdot (\nabla\times B) dx\\
&\leq \frac12\int\left|\nabla\times B\right|^2 dx+C(\left|u\right|_\infty)\int \left|B\right|^2 dx.
\end{aligned}
\end{equation}
Due to the assumption 3, we know that $\nu_b(\rho)$ has  lower bound, then by Gronwall inequality to \eqref{e3.9} we can get that $B=0$ .Furthermore, there exists a continuous solution operator $G$ from $C([0,T];X_n)$ to $L^\infty(0,T;L^2)\cap L^2(0,T;H^1)$ by $G(u)=B$.

Now we turn to solve the approximate momentum equation on the space $X_n$. For $\rho=S(u), B=G(u)$ , we are looking for a function $u\in C([0,T];X_n)$ such that 
\begin{equation}\label{e3.10}
\begin{aligned}
&\int_{\O} \rho u(T)\varphi dx-\int_{\O}m_0\varphi dx-\int_0^T\int_{\O}\left(\rho u\otimes u:\nabla\varphi-P(\rho)\diver\varphi-P_c(\rho)\diver\varphi\right) dxdt\\
&+2\kappa^2\int_0^T\int_{\O}\frac{\Delta\sqrt\rho}{\sqrt\rho}\diver(\rho\varphi)dxdt+\eps\int_0^T\int_{\O}\nabla\rho\cdot\nabla u\cdot\varphi dxdt\\
&-\delta\int_0^T\int_{\O}\Delta^s(\diver(\rho\varphi)):\Delta^{s+1}\rho dxdt+2\int_0^T\int_{\O}\rho D(u)\cdot\nabla\varphi dxdt\\
&+\eta\int_0^T\int_{\O}\Delta u\cdot\Delta\varphi dxdt -\int_0^T\int_{\O}\left((\nabla\times B)\times B\right)\cdot\varphi dxdt=0,
\end{aligned}
\end{equation}
for all $\varphi\in X_n$.
we will apply the Banach fixed point theorem to prove the local existence of solutions for the equation \eqref{e3.10}. Following the same argument \cite{Feireisl2001} we can solve \eqref{e3.10}. Next, we introduce an operator defined on the set $\{\rho\in L^1, \rho\ge \underline\rho>0\}$, and $M[\rho]:X_n\to X_n^{\ast} $
\begin{equation}
 <M[\rho]v,u>=\int \rho v\cdot u dx , v,u\in X_n,
\end{equation}
where $X_n^{\ast}$ stands for the dual space of $X_n$. And the operator $M[\rho]$ has the following properties:
\begin{itemize}
\item $\left\|M[\rho]\right\|_{L(X_n,X_n^{\ast})}\leq C(n)\left\|\rho\right\|_{L^1}$.
\item $M^{-1}[\rho]$ is invertible under the condition $\rho\ge \underline\rho>0$ and 
 $$\left\|M^{-1}[\rho]\right\|_{L(X_n^{\ast},X_n)}\leq (\underline\rho)^{-1},$$
where $M^{-1}[\rho]:X_n^{\ast}\to X_n$.
\item $M^{-1}[\rho]$ is Lipschtiz continuous.
\end{itemize}
\end{enumerate}
The first two properties of the operator is easily to get. Since$$M^{-1}[\rho_1]-M^{-1}[\rho_2]=M^{-1}[\rho_2](M[\rho_2]-M[\rho_1])M^{-1}[\rho_1],$$
where $\rho_1,\rho_2\in\{\rho\in L^1, \rho\ge\underline\rho>0\}$. we can get 
$$\left\|M^{-1}[\rho_1]-M^{-1}[\rho_2]\right\|_{L(X_n^{\ast},X_n)}\leq C(n,\underline\rho)\left\|\rho_1-\rho_2\right\|_{L^2}.$$
Thus, $M^{-1}$ is Lipschitz continuous.
Now, using the definition of the operator $M[\rho]$ the equation \eqref{e3.10} can be rewritten as 
\begin{equation}
u_n(t)=M^{-1}[\rho]\left(m_0^{\ast}+\int_0^T N[\rho(s),u(s),B(s)] ds\right).
\end{equation}
where $\rho=S(u), B=G(u)$ and 
\begin{equation}
\begin{aligned}
 N[\rho(s),u(s),B(s)] &=\int_{\O}\left(\rho u\cdot\varphi_t+\rho u\otimes u:\nabla\varphi+P(\rho)\diver\varphi+P_c(\rho)\diver\varphi\right)dx\\
&-2\kappa^2\int_{\O}\frac{\Delta\sqrt\rho}{\sqrt\rho}\diver(\rho\varphi)dx-\eps\int_{\O}\nabla\rho\cdot\nabla u\cdot\varphi dx\\
&-\delta\int_{\O}\Delta^s(\diver(\rho\varphi)):\Delta^{s+1}\rho dx-2\int_{\O}\rho D(u)\cdot\nabla\varphi dx\\
&-\eta\int_{\O}\Delta u\cdot\Delta\varphi dx+\int_{\O}\left((\nabla\times B)\times B\right)\cdot\varphi dx, \varphi\in X_n.
\end{aligned}
\end{equation}
Due to the operators $S,G,M^{-1}$ are Lipschitz continuous, the above nonlinear equation can be solved by the fixed point argument on a short time interval $[0,T^{'}], T^{'}\leq T$ in the Banach space $C([0,T];X_n)$. Therefore, we can proved the local existence of solutions $(\rho_n,u_n,B_n)$ to the approximate systems \eqref{e3.3},\eqref{e3.8} and \eqref{e3.10}.
\subsection{Uniform estimate and global existence of solutions}

Assume $(\rho_n, u_n,B_n)$ is the approximate solutions exists on the $[0,T^{'}], T^{'}\leq T$. Our goal in this section is to extend the approximate solutions $(\rho_n, u_n,B_n)$ to the whole interval $[0,T]$, it is sufficient to establish the uniform bound on the norm $\left\|u_n\right\|_{X_n}$ which allows us to iterate the above procedure many times to reach the whole interval $[0,T]$.
\begin{lem}
Assume $T^{'}\leq T$, $(\rho_n,u_n,B_n)$ be the solutions to the \eqref{e3.3},\eqref{e3.8} and \eqref{e3.10}.  then we have the following holds:
\begin{equation}\label{e3.14}
\begin{aligned}
&\frac{d}{dt}E(\rho_n,u_n,B_n)+2\int_{\O} \rho_n \left|D(u_n)\right|^2 dx +\eps\int_{\O}\left(H^{''}(\rho_n)+H^{''}_c(\rho_n)\right)\left|\nabla\rho_n\right|^2 dx\\
&+\int_{\O}\nu_b(\rho_n)\left|\nabla\times B_n\right|^2 dx +\eta\int_{\O}\left|\Delta u_n\right|^2 dx+\delta\eps\int_{\O}\left|\Delta^{s+1}\rho_n\right|^2 dx\\
&+\eps\kappa^2\int_{\O} \rho_n\left|\nabla^2\log\rho_n\right|^2 dx=0,
\end{aligned}
\end{equation}
where 
\begin{equation}
E(\rho_n,u_n,B_n)=\int_{\O}\left(\frac12\rho_n\left|u_n\right|^2+H(\rho_n)+H_c(\rho_n)+\kappa^2\left|\nabla\sqrt\rho_n\right|^2+\frac12\left|B_n\right|^2+\frac{\delta}{2}\left|\nabla^{2s+1}\rho_n\right|^2\right)dx.
\end{equation}
\end{lem}
\begin{pf}
 Differentiating \eqref{e3.10} with respect to time and using the test function $\varphi=u_n$ we have 
 \begin{equation}\label{e3.16}
 \begin{aligned}
 &\frac{d}{dt}\int \frac12\rho\left|u\right|^2dx -\int(\d\rho+ \diver(\rho u))\frac{\left|u\right|^2}{2}dx+\eps\int \Delta\rho \left|u\right|^2 dx+\eps\int \nabla\rho\cdot\nabla u\cdot u dx\\
 &+\int\nabla(P(\rho)+P_c(\rho))\cdot u dx-\delta\int_0^T\int_{\O}\Delta^s(\diver(\rho u)):\Delta^{s+1}\rho dxdt+2\int_0^T\int_{\O}\rho \left|D(u)\right|^2 dxdt\\
 &\eta\int\left|\Delta u\right|^2 dx+2\kappa^2\int \frac{\Delta\sqrt\rho}{\sqrt\rho}\diver(\rho u)dxdt-\int \left(\nabla\times B\right)\times B\cdot u dx=0.
 \end{aligned}
 \end{equation}
 due to the fact that 
 \begin{equation}\label{e3.17}
 \begin{aligned}
 &\int\nabla(P(\rho)+P_c(\rho))\cdot u dx\\
 &=\int\frac{1}{\rho}(P^{'}(\rho)+P_c^{'}(\rho))\nabla\rho\cdot \rho u dx\\
 &=\int\nabla(H^{'}(\rho)+H_c^{'}(\rho))\cdot \rho u dx\\
 &=-\int (H^{'}(\rho)+H_c^{'}(\rho))\diver(\rho u) dx\\
 &=\int  (H^{'}(\rho)+H_c^{'}(\rho))(\d\rho-\eps\Delta\rho)dx\\
 &=\frac{d}{dt}\int(H(\rho)+H_c(\rho)) dx+\eps\int (H^{''}(\rho)+H_c^{''}(\rho))\left|\nabla\rho\right|^2 dx,
 \end{aligned}
 \end{equation}
 \begin{equation}\label{e3.18}
 \begin{aligned}
 &2\kappa^2\int\frac{\Delta\sqrt\rho}{\sqrt\rho}\diver(\rho u)dxdt\\
 &=2\kappa^2\int\frac{\Delta\sqrt\rho}{\sqrt\rho}\Delta\rho-4\kappa^2\int\Delta\sqrt\rho\d\sqrt\rho dx\\
 &=\eps\kappa^2\int\rho\left|\nabla^2\log\rho\right|^2 dx +\frac{\kappa^2}{2}\frac{d}{dt}\int\left|\nabla\sqrt\rho\right|^2 dx.
 \end{aligned}
 \end{equation}
and similarly calculation
\begin{equation}\label{e3.19}
\delta\int_{\O}\Delta^s(\diver(\rho u)):\Delta^{s+1}\rho dx=\eps\delta\int\left|\Delta^{s+1}\rho\right|^2 dx+\frac{\delta}{2}\frac{d}{dt}\int \left|\nabla^{2s+1}\rho\right|^2 dx.
\end{equation} 
Then multiply the magnetic equation by $B_n$ we get
\begin{equation}\label{e3.20}
\frac12\int \left|B_n\right|^2-\int\nabla\times(u\times B_n)\cdot B_n dx+\int \nu_b(\rho)\left|\nabla\times B_n\right|^2 dx=0.
\end{equation}
where we used the following the identity
\begin{equation*}
\int (\nabla\times B)\times B\cdot u dx=-\int\nabla\times(u\times B)\cdot B dx.
\end{equation*}
Summing up \eqref{e3.16} and \eqref{e3.20} together, we can get the desire estimate \eqref{e3.14}.
\end{pf}
From Lemma 3.1 we have the following estimate
\begin{equation}\label{e3.21}
\begin{aligned}
&\rho_n\in L^\infty(0,T;L^\gamma),\rho_n^{-1}\in L^\infty(0,T;L^{\gamma^{-}}),\nabla\sqrt\rho_n\in L^\infty(0,T;L^2),\\
&\sqrt\rho_n u_n\in L^\infty(0,T;L^2),\sqrt\rho_n Du_n\in L^2(0,T;L^2),\sqrt\eta\Delta u_n\in L^2(0,T;L^2),\\
&\nabla\rho_n^{\frac{\gamma}{2}}\in L^2(0,T;L^2),\sqrt\delta \rho_n\in L^\infty(0,T;H^{2s+1}),\\
&B_n\in L^\infty(0,T;L^2) \cap L^2(0,T;H^1),\sqrt{\delta\eps}\rho_n\in L^2(0,T;H^{2s+2}),\\
&\sqrt\eps\sqrt\rho_n\nabla^2\log\rho_n\in L^2(0,T;L^2),\sqrt{\nu_b(\rho_n)}\nabla B_n\in L^2(0,T;L^2).
\end{aligned}
\end{equation}
By sobolev embedding we have $\left\|\rho^{-1}\right\|_{L^\infty}\leq C\left\|\rho^{-1}\right\|_{W^{3,1}}$ and 
\begin{equation}
\left\|\nabla^3\rho^{-1}\right\|_{L^1}\leq C\left(1+\left\|\nabla^3\rho\right\|_{L^\infty(0,T;L^2)}\right)^3\left(1+\left\|\rho^{-1}\right\|_{L^\infty(0,T;L^4)}\right)^4.
\end{equation}
Therefore, provided that $\gamma^{-}\ge 4, 2s+1\ge 3$ we can get 
\begin{equation}\label{e3.23}
\left\|\rho^{-1}\right\|_{L^\infty((0,T)\times \O)} \leq C(\delta).
\end{equation}
Furthermore, by lemma 2.4 we get 
\begin{equation}\label{e3.24}
\sqrt\eps \kappa\left\|\sqrt\rho_n\right\|_{L^2(0,T;H^2)}+\eps^{1/4}\kappa^{1/2}\left\|\nabla\rho_n^{1/4}\right\|_{L^4(0,T;L^4)}\leq C.
\end{equation}
Together with \eqref{e3.21} we get the uniform bound for $u_n$, i.e
$$\sup_{[0,T_{\max}]}\left\|u_n\right\|_{X_n}\leq C.$$ where $C$ is independent of $T_{\max}$. Thus, we get a global existence of approximate solutions.
\section{passage to the limit $n\to\infty$}
After we have already constructed a family of approximate solutions $(\rho_n,u_n,B_n)$. The purpose of this section is to let $n\to\infty$. This can be achieved by using the uniform estimate of approximate solutions and Aubin-Lions lemma.
\begin{lem}
The following estimates holds for some constant $C$ independent of $n$
\begin{equation}\label{e4.1}
\begin{aligned}
&\left\|\d\rho_n\right\|_{L^2(0,T;L^2)}+\left\|\rho_n\right\|_{L^2(0,T;H^{2s+2})}\leq C,\\
&\left\|\d\sqrt\rho_n\right\|_{L^2(0,T;H^{-1})}+\left\|\sqrt\rho_n\right\|_{L^2(0,T;H^{2})}\leq C,\\
&\left\|\d(\rho_n u_n)\right\|_{L^2(0,T;H^{-(2s+1)})}+\left\|(\rho_n u_n)\right\|_{L^2(0,T;W^{1,3/2})}\leq C,\\
&\left\|P(\rho_n)\right\|_{L^{5/3}}+\left\|P_c(\rho_n)\right\|_{L^{5/3}}\leq C,\\
&\left\|\d B_n\right\|_{L^2(0,T;H^{-1})}+\left\|B_n\right\|_{L^2(0,T;H^1)}\leq C,\\
&\left\|\d(1/\sqrt\rho_n)\right\|_{L^\infty(0,T;W^{-1,1})}\leq C.
\end{aligned}
\end{equation}
\end{lem}
\begin{pf}
By the continuity equation we have 
\begin{equation*}
\d\rho_n=\eps\Delta\rho_n-\diver(\rho_n u_n)\in L^2(0,T;L^2).
\end{equation*}
Then since 
\begin{equation*}
\d\sqrt\rho_n+\frac{1}{2\sqrt\rho_n}\diver(\rho_n u_n)=\eps\left(\Delta\sqrt\rho_n+\frac{\left|\nabla\sqrt\rho_n\right|^2}{\sqrt\rho_n}\right).
\end{equation*}
together with the estimate \eqref{e3.24}  yields that 
$$\d\sqrt\rho_n\in L^2(0,T;H^{-1}).$$
due to the momentum equation
\begin{align*}
\d(\rho_n u_n)&=-\diver(\rho_n u_n\otimes u_n)-\nabla(P(\rho_n)+P_c(\rho_n))+2\diver(\rho_n D(u_n))+\eta\Delta u_n\\
&\quad +2\kappa^2\rho_n\nabla\left(\frac{\Delta\sqrt\rho_n}{\sqrt\rho_n}\right)-(\nabla\times B_n)\times B_n+\delta\rho_n\nabla\Delta^{2s+1}\rho_n+\eps\nabla\rho_n\cdot\nabla u_n.
\end{align*}
we deduce that $\d(\rho_n u_n)\in L^2(0,T;H^{-(2s+1)})$.
Moreover,  we also have
$$\nabla(\rho_n u_n)=2\nabla\sqrt\rho_n\sqrt\rho_n u_n+2\sqrt\rho_n\nabla u_n\nabla\sqrt\rho_n\in L^2(0,T;L^{3/2}).$$
then we have $P(\rho_n)\in L^\infty(0,T;L^1)$ and $P(\rho_n)\in L^1(0,T;L^3)$ by the interpolation inequality we get 
$P(\rho_n)\in L^{5/3}(0,T;L^{5/3})$. and similarly we can proof $P_c(\rho_n)$.Finally, by the magnetic equation and \eqref{e3.24} we have 
$$\d B=\nabla\times(u\times B)-\nabla\times(\nu_b(\rho)\nabla\times B)\in L^2(0,T;H^{-1}).$$
Since
 $$\d\left(\frac{1}{\sqrt\rho_n}\right)+\frac12\diver(\rho^{-1/2}u)-\frac32\rho^{-3/2}\diver u=-\eps\left(\frac{\Delta\sqrt\rho_n}{\sqrt\rho_n}+\frac{\left|\nabla\sqrt\rho_n\right|^2}{\rho_n^{3/2}}\right).$$
 and using the previous estimate \eqref{e3.21} we have  $\d\rho_n^{-1/2}\in L^2(0,T;W^{-1,6/5})$.
\end{pf}
By Aubin-Lions lemma and the Lemma 4.1, we have 
\begin{equation}\label{e4.2}
\begin{aligned}
&\sqrt\rho_n\to \sqrt\rho \,\text{ strongly in} \,L^2(0,T;H^1),\\
&\rho_n\to\rho\,\text{strongly in } \,L^2(0,T;H^{2s+1}), \text{and weakly in } L^2(0,T;H^{2s+2}),\\
&\rho_n u_n\to\rho u \,\text{strongly in}\, L^2(0,T;L^2),\\
&u_n\to u\,\text{weakly in}\, L^2(0,T;L^2), B_n\to B\, \text{strongly in}\, L^2(0,T;L^2),\\
&\rho_n^{-1/2}\to \rho^{-1/2} almost \,everywhere.\nabla B_n\to \nabla B\, \text{weakly in}\, L^2(0,T;L^2).
\end{aligned}
\end{equation}
Moreover, we can get that $ \rho_n u_n\otimes u_n\to \rho u\otimes u$ in the distribution sense. Together with the estimate \eqref{e3.21} and \eqref{e4.2}, we can prove that 
\begin{align*}
&P(\rho_n)\to P(\rho) \,\text{strongly in }\, L^1((0,T)\times\O),\\
&P_c(\rho_n)\to P_c(\rho) \,\text{strongly in }\, L^1((0,T)\times\O).
\end{align*}
The viscosity term can pass to the limit 
$$\iint\diver(\rho_n D(u_n) \varphi dxdt\to \iint\diver(\rho D(u) \varphi dxdt.$$
due to the fact that
\begin{equation*}
\begin{aligned}
&\iint\rho(\nabla u+(\nabla u)^{\top}):\nabla\varphi dxdt\\
&=\iint(\rho\partial_i u^j\partial_i\varphi_j +\rho\partial_j u^i\partial_i\varphi_j)dx dt\\
&=\iint(\rho u^j)_i\partial_i\varphi_j+(\rho u^i)_j\partial_i\varphi^j -\partial_i\rho u^j\partial_i\varphi^j-\partial_j\rho u^i\partial_i\varphi^j dx\\
&=-2\iint(\nabla\sqrt\rho\otimes\sqrt\rho u)\cdot\nabla\varphi dxdt-2\iint(\sqrt\rho u\otimes\nabla\sqrt\rho):\nabla\varphi dxdt\\
&-\iint\sqrt\rho\sqrt\rho u\cdot \Delta\varphi dxdt-\iint\sqrt\rho\sqrt\rho u\nabla\diver\varphi dxdt.
\end{aligned}
\end{equation*}
For the quantum term can also pass to the limit by using the convergence of the $\sqrt\rho_n$. In fact, for any test function
\begin{equation}\label{e4.3}
\begin{aligned}
&\int \rho_n\nabla\left(\frac{\Delta\sqrt\rho_n}{\sqrt\rho_n}\right) \varphi dx\\
&=-\int \frac{\Delta\sqrt\rho_n}{\sqrt\rho_n}\diver(\rho_n\varphi) dx\\
&=-\int \Delta\sqrt\rho_n\sqrt\rho_n\diver\varphi dx-2\int\Delta\sqrt\rho_n\nabla\sqrt\rho_n\varphi dx\\
&=-\int\diver(\sqrt\rho_n\nabla\sqrt\rho_n)\diver\varphi dx+\int\left|\nabla\sqrt\rho_n\right|^2\diver\varphi dx\\
&\quad-2\int \diver(\nabla\sqrt\rho_n\otimes\nabla\sqrt\rho_n)\varphi dx+2\int(\nabla\sqrt\rho_n\cdot\nabla)\nabla\sqrt\rho_n\varphi dx\\
&=\int \sqrt\rho_n\nabla\sqrt\rho_n\nabla\diver\varphi dx+2\int \nabla\sqrt\rho_n\otimes\nabla\sqrt\rho_n\nabla\varphi dx.
\end{aligned}
\end{equation}
For the capillarity term  we can pass to the limit $$\delta\iint\rho_n\nabla\Delta^{2s+1}\rho_n\varphi dxdt\to  \delta \iint \rho \nabla \Delta^{2s+1} \rho\varphi dxdt.$$
due to the strong convergence of $\rho_n$ and the following fact
\begin{equation*}
\delta\iint\rho_n\nabla\Delta^{2s+1}\rho_n\varphi dxdt=\delta\iint\nabla^{2s+1}\rho_n\Delta^{s}\diver(\rho_n\varphi) dxdt.
\end{equation*}
 From the assumption, we get $\nu_b(\rho)$ has uniform lower bound which yields the weak convergence of $\nabla B_n$. Then together with the strong convergence of the $B_n$ and weak convergence of $u_n$ enable us pass to the limit for the magnetic equation. Therefore, the proof of pass to limit $n\to\infty$ is completed.Thus, we can show that $(\rho,u,B)$ solve the following systems
 \begin{equation}\label{e4.4}
\left\{
\begin{aligned}
&\d \rho+\diver(\rho u)=\eps\Delta\rho, x\in \O, t>0,\\
&\d(\rho u)+\diver(\rho u\otimes u)+\nabla(P(\rho)+P_c(\rho))-2\diver(\rho D(u))\\
&\quad+\eta\Delta^2 u+\eps\nabla\rho\cdot \nabla u-\delta\rho\nabla\Delta^{2s+1}\rho -2\kappa^2\rho\nabla\left(\frac{\Delta\sqrt\rho}{\sqrt\rho}\right)\\
&\quad-(\nabla\times B)\times B=0,\\
&\d B-\nabla\times(u\times B)+\nabla(\nu_b(\rho)\nabla\times B)=0.\\
\end{aligned}
\right.
\end{equation}
 in the distribution sense and also satisfies the energy estimate \eqref{e3.14} and \eqref{e3.24}.
 \section{B-D entropy estimate and pass to the limit $\eps,\eta,\delta\to 0$}
 The purpose of this section is to derive the B-D entropy estimate for the approximated system \eqref{e4.4}. This estimate first established by Bresch-Desjardin-Lin in \cite{Didier2003On}. By the \eqref{e3.23} and \eqref{e4.1} we get 
 \begin{equation}
 \rho(x,t)\ge C(\delta)>0,\text{ and}\, \rho\in L^2(0,T;H^{2s+2})\cap L^\infty(0,T;H^{2s+1}).
 \end{equation}
then we can choose $\nabla\phi(\rho)=2\nabla\rho/\rho$ as test function to multiply the momentum equation to derive the entropy estimate. 
\begin{lem}(B-D entropy estimate)
The following equality holds 
\begin{equation}\label{e5.2}
\begin{aligned}
&\frac{d}{dt}\int\left(\frac12\rho\left|u+\nabla\phi(\rho)\right|^2+H(\rho)+H_c(\rho)+\kappa^2\left|\nabla\sqrt\rho\right|^2+\frac12\left|B\right|^2+\frac{\delta}{2}\left|\nabla^{2s+1}\rho\right|^2\right)dx\\
&\quad+\eta\int \left|\Delta u\right|^2 dx+2\int \rho\left|A(u)\right|^2 dx+2\int\frac{1}{\rho}(P^{'}(\rho)+P_c^{'}(\rho))\left|\nabla\rho\right|^2 dx\\
&\quad+2\kappa^2\int\rho\left|\nabla^2\log\rho\right|^2 dx+\eps\kappa^2\int\rho\left|\nabla^2\log\rho\right|^2 dx+\int\nu_b(\rho)\left|\nabla\times B\right|^2dx\\
&\quad+\eps\delta\int\left|\Delta^{s+1}\rho\right|^2 dx+2\delta\int\rho\left|\Delta^{s+1}\rho\right|^2dx+\eps\int\frac{1}{\rho}(P^{'}(\rho)+P_c^{'}(\rho))\left|\nabla\rho\right|^2 dx\\
&=\eps\int \nabla\phi(\rho)\cdot\nabla(\phi^{'}(\rho)\Delta\rho) dx-\eps\int \nabla\rho\cdot\nabla u\cdot\nabla\phi(\rho) dx+\eps\int \frac{\left|\nabla\phi(\rho)\right|^2}{2}\Delta\rho dx\\
&\quad-\eta\int \Delta u\cdot\nabla\Delta\phi(\rho) dx-\eps\int \diver(\rho u)\phi^{'}(\rho)\Delta\rho dx+\int(\nabla\times B)\times B\cdot \nabla\phi(\rho) dx.
\end{aligned}
\end{equation}
\end{lem}
\begin{pf}
We first multiply the approximate continuity equation by $\frac{\left|\nabla\phi(\rho)\right|^2}{2}$, we have
\begin{equation}
\begin{aligned}
&\frac{d}{dt}\int\frac12\rho\left|\nabla\phi(\rho)\right|^2 dx\\
&=\int\rho\d\left(\frac{\left|\nabla\phi(\rho)\right|^2}{2}\right) dx+\int\d\rho\frac{\left|\nabla\phi(\rho)\right|^2}{2} dx\\
&=\int\rho\nabla\phi(\rho)\cdot\nabla(\phi^{'}(\rho)\d\rho) dx+\int\d\rho\frac{\left|\nabla\phi(\rho)\right|^2}{2}dx\\
&=\int\rho\nabla\phi(\rho)\cdot\nabla(\phi^{'}(\rho)\d\rho) dx+\eps\int \Delta\rho\frac{\left|\nabla\phi(\rho)\right|^2}{2} dx-\int\frac{\left|\nabla\phi(\rho)\right|^2}{2} \diver(\rho u) dx.\\
\end{aligned}
\end{equation}
due to the first term on the right hand side is equal to 
\begin{equation}
\begin{aligned}
&\int\rho\nabla\phi(\rho)\cdot\nabla(\phi^{'}(\rho)\d\rho) dx\\
&=\eps\int\rho\nabla\phi(\rho)\nabla(\phi^{'}(\rho)\Delta\rho) dx-\int\rho\nabla u:\nabla\phi(\rho)\otimes \nabla\phi(\rho)dx+\int\rho\left|\nabla\phi(\rho)\right|^2\diver u dx\\
&\quad+\int\rho^2\phi^{'}(\rho)\Delta\phi(\rho)\diver u dx+\int \left|\nabla\phi(\rho)\right|^2\diver(\rho u) dx+\int\rho u\cdot\nabla^2\phi(\rho)\cdot\nabla\phi(\rho) dx\\
&=\eps\int\rho\nabla\phi(\rho)\nabla(\phi^{'}(\rho)\Delta\rho) dx-\int\rho\nabla u:\nabla\phi(\rho)\otimes \nabla\phi(\rho)dx+\int\rho\left|\nabla\phi(\rho)\right|^2\diver u dx\\
&\quad+\int\rho^2\phi^{'}(\rho)\Delta\phi(\rho) \diver u dx+\frac12\int\left|\nabla\phi(\rho)\right|^2\diver(\rho u) dx.
\end{aligned}
\end{equation}
Thus 
\begin{equation}\label{e5.5}
\begin{aligned}
&\frac{d}{dt}\int\frac12\rho\left|\nabla\phi(\rho)\right|^2 dx=\eps\int\rho\nabla\phi(\rho)\nabla(\phi^{'}(\rho)\Delta\rho) dx-\int\rho\nabla u:\nabla\phi(\rho)\otimes \nabla\phi(\rho)dx\\
&\quad+\int\rho\left|\nabla\phi(\rho)\right|^2\diver u dx+\int\rho^2\phi^{'}(\rho)\Delta\phi(\rho) \diver u dx+\eps\int \Delta\rho\frac{\left|\nabla\phi(\rho)\right|^2}{2} dx.
\end{aligned}
\end{equation}
Next, we need to calculation the following term
\begin{equation}\label{e5.6}
\begin{aligned}
&\frac{d}{dt}\int \rho u\cdot\nabla\phi(\rho) dx\\
&=\int \d(\rho u)\cdot \nabla\phi(\rho) dx+\int \rho u\cdot\nabla(\phi^{'}(\rho)\d\rho) dx\\
&=\int \d(\rho u)\cdot \nabla\phi(\rho) dx+\eps\int\rho u\cdot\nabla(\phi^{'}(\rho)\Delta\rho) dx+\int\phi^{'}(\rho)(\diver(\rho u))^2 dx.
\end{aligned}
\end{equation}
Combine \eqref{e5.5} and \eqref{e5.6} together, we have 
\begin{equation}\label{e5.7}
\begin{aligned}
&\frac{d}{dt}\int\frac12\rho\left|\nabla\phi(\rho)\right|^2+\int \rho u\cdot\nabla\phi(\rho) )dx\\
&=\eps\int\rho\nabla\phi(\rho)\nabla(\phi^{'}(\rho)\Delta\rho) dx-\int\rho\nabla u:\nabla\phi(\rho)\otimes \nabla\phi(\rho)dx+\int\rho\left|\nabla\phi(\rho)\right|^2\diver u dx\\
&\quad+\int\rho^2\phi^{'}(\rho)\Delta\phi(\rho) \diver u dx+\eps\int \Delta\rho\frac{\left|\nabla\phi(\rho)\right|^2}{2}+\int\phi^{'}(\rho)(\diver(\rho u))^2 dx\\
&\quad+\int\diver(\rho u\otimes u):\nabla\phi(\rho) dx-\eps\int\diver(\rho u)\phi^{'}(\rho)\Delta\rho dx-\int\nabla(P(\rho)+P_c(\rho))\cdot\nabla\phi(\rho) dx\\
&\quad+2\int\nabla u:\nabla\rho\otimes\nabla\phi(\rho) dx-2\int\nabla\rho\cdot\phi(\rho)\diver u dx-\eps\int\nabla\rho\cdot\nabla u\cdot\nabla\phi(\rho) dx\\
&\quad-2\int\rho\Delta\phi(\rho)\diver u dx-2\kappa^2\int\rho\left|\nabla^2\log\rho\right|^2 dx-2\delta\int\left|\Delta^{2s+1}\rho\right|^2 dx\\
&\quad-\eta\int\Delta u\cdot\nabla\Delta\phi(\rho) dx+\int((\nabla\times B)\times B)\cdot\nabla\phi(\rho) dx.\\
&=\eps\int\rho\nabla\phi(\rho)\nabla(\phi^{'}(\rho)\Delta\rho) dx+\eps\int \Delta\rho\frac{\left|\nabla\phi(\rho)\right|^2}{2}-\int\nabla(P(\rho)+P_c(\rho))\cdot\nabla\phi(\rho) dx\\
&\quad+2\int\rho\left|D(u)\right|^2 dx-2\int \rho \left|A(u)\right|^2 dx-\eps\int\diver(\rho u)\phi^{'}(\rho)\Delta\rho dx-2\delta\int\left|\Delta^{2s+1}\rho\right|^2 dx\\
&\quad-2\kappa^2\int\rho\left|\nabla^2\log\rho\right|^2 dx-\eta\int\Delta u\cdot\nabla\Delta\phi(\rho) dx-\eps\int\nabla\rho\cdot\nabla u\cdot\nabla\phi(\rho) dx\\
&\quad+\int((\nabla\times B)\times B)\cdot\nabla\phi(\rho) dx.
\end{aligned}
\end{equation}
Together with \eqref{e5.7} and \eqref{e3.14} yields the desire estimate. Therefore,the proof of lemma is completed.
\end{pf}
Next, we turn to estimate each term on the right hand side of \eqref{e5.2}. The sign of the first term is negative. In fact, integration by parts 
\begin{equation}
\eps\int \nabla\phi(\rho)\cdot\nabla(\phi^{'}(\rho)\Delta\rho) dx=-4\eps\int \frac{1}{\rho}\left|\Delta\rho\right|^2 dx.
\end{equation}
For the second term which can be estimate as follows
\begin{align}
\left|\eps\int \nabla\rho\cdot\nabla u\cdot\nabla\phi(\rho) dx\right|&=2\eps\left|\int\frac{1}{\rho}\nabla\rho\cdot\nabla u\cdot\nabla\rho dx\right|\\
&=2\eps\left|\int \frac{\left|\nabla\rho\right|^2}{\rho^{3/2}}\rho^{1/2}\nabla u dx\right|\notag\\
&\leq \eps\int \frac{\left|\nabla\rho\right|^4}{\rho^3} dx+\eps\int\rho\left|\nabla u\right|^2 dx\notag.
\end{align}
due to 
\begin{equation*}
\eps\int\rho\left|\nabla u\right|^2 dx\leq 2\eps(\int \rho\left|A(u)\right|^2 dx+\int\rho\left|D(u)\right|^2 dx),
\end{equation*}
Therefore, we have 
\begin{equation}
\left|\eps\int \nabla\rho\cdot\nabla u\cdot\nabla\phi(\rho) dx\right|\leq \eps\int \frac{\left|\nabla\rho\right|^4}{\rho^3} dx+2\eps\int \rho\left|A(u)\right|^2 dx+2\eps\int\rho\left|D(u)\right|^2 dx.
\end{equation}
The third term can be estimated as follows
\begin{align}
\left|\eps\int \frac{\left|\nabla\phi(\rho)\right|^2}{2}\Delta\rho dx\right|&=2\eps\left|\frac{\left|\nabla\rho\right|^2}{\rho^2}\Delta\rho dx\right|\\
&=2\eps\left|\int\frac{\left|\nabla\rho\right|^2}{\rho^{3/2}}\frac{\Delta\rho}{\rho^{1/2}} dx\right|\notag\\
&\leq \eps\int \frac{\left|\nabla\rho\right|^4}{\rho^3}dx +\eps\int \frac{\left|\Delta\rho\right|^2}{\rho} dx\notag.
\end{align}
The fourth term estimated as follows
\begin{equation}
\left|\eta\int \Delta u\cdot\nabla\Delta\phi(\rho) dx\right|\leq \frac{\eta}{2}\left\|\Delta u\right\|_2^2+\frac{\eta}{2}\left\|\nabla\Delta\phi(\rho)\right\|_2^2.
\end{equation}
and the second term on the right hand side is equal to 
\begin{align*}
\nabla\Delta\phi(\rho)&=2\partial_{kk}\left(\frac{\partial_i \rho}{\rho}\right)\\
&=2\partial_{k}\left(-\frac{1}{\rho^2}\partial_{k}\rho\partial_{i}\rho+\frac{\partial_{ik}\rho}{\rho}\right)\\
&=2\left(2\rho^{-3}(\partial_{k}\rho)^2\partial_i\rho-\frac{1}{\rho^2}\partial_{kk}\rho\partial_i\rho-\frac{2}{\rho^2}\partial_k\rho\partial_{ik}\rho+\frac{\partial_{ikk}\rho}{\rho}\right)\\
&=\frac{2\nabla\Delta\rho}{\rho}-\frac{4(\nabla\rho\cdot\nabla)\nabla\rho}{\rho^2}+\frac{4\left|\nabla\rho\right|^2\Delta\rho}{\rho^3}-\frac{2\Delta\rho\nabla\rho}{\rho^2}.
\end{align*}
then we get
\begin{equation}
\left\|\nabla\Delta\phi(\rho)\right\|_2\leq C(1+\left\|\rho\right\|_{H^{2s+1}})^3(1+\left\|\rho^{-1}\right\|_{L^\infty})^3\leq C(\delta).
\end{equation}
we can choose $\eta$ small enough with respect to $\delta$ such that we can get the uniform bound.

The fifth term on the right hand side of \eqref{e5.2} is equal to 
\begin{equation}
2\eps\int \diver(\rho u)\phi^{'}(\rho)\Delta\rho dx =2\eps\int \Delta\rho\diver u dx+2\eps\int\frac{\Delta\rho}{\rho}\nabla\rho\cdot udx.
\end{equation} 
since 
\begin{align*}
2\eps\left|\int \Delta\rho\diver u dx\right|&\leq \eps\int\frac{\left|\Delta\rho\right|^2}{\rho} dx+\eps\int \rho\left|\nabla u\right|^2 dx\\
&\leq \eps\int\frac{\left|\Delta\rho\right|^2}{\rho} dx+2\eps\int \rho\left|A(u)\right|^2 dx+\int\rho\left|D(u)\right|^2 dx.\\
\end{align*}
\begin{equation*}
2\eps\left|\int\frac{\Delta\rho}{\rho}\nabla\rho\cdot udx\right|\leq2\eps\left\|\rho^{-1}\right\|_{L^\infty}^{3/2}\left\|\rho\right\|_{H^{2s+1}}.
\end{equation*}
we can choose $\eps$ small enough with respect to $\delta$ yields the uniform bound.

Finally, we estimate the last term on the right hand side of \eqref{e5.2}.
\begin{equation}
\int(\nabla\times B)\times B\cdot \nabla\phi(\rho) dx\leq \int \frac{\left|\nabla\times B\right|^2}{\eps\rho^2} dx+\eps\int \left|\nabla\rho\times B\right|^2 dx.
\end{equation}
The first term can be absorbed by the magnetic diffusion term under the Assumption 3, and the second term estimated as follows
\begin{equation*}
\eps\int\left|\nabla\rho\times B\right|^2 dx\leq \eps\left\|\nabla\rho\right\|_{L^\infty}^2\left\|B\right\|_2^2\leq \left\|\nabla\rho\right\|_{H^2}^2\left\|B\right\|_2^2\leq C\eps\left\|\rho\right\|_{H^{2s+1}}^2.
\end{equation*}
 we can follow the same argument choose $\eps$ small enough with respect to $\delta$ enable us get the uniform bound.
 \subsection{pass to the limits as $\eps,\eta\to 0$}
 
 In this step we assume $(\rho_{\eps,\eta},u_{\eps,\eta},B_{\eps,\eta})$ are the approximate solutions.  Then from the B-D entropy estimate we can deduce that
 \begin{equation}\label{e5.16}
 \begin{aligned}
 &\rho_{\eps,\eta}\in L^\infty(0,T;L^\gamma),\rho_{\eps,\eta}^{-1}\in L^\infty(0,T;L^{\gamma^{-}}),\nabla\sqrt\rho_{\eps,\eta}\in L^\infty(0,T;L^2),\\
&\sqrt\rho_{\eps,\eta} u_{\eps,\eta}\in L^\infty(0,T;L^2),\sqrt\rho_{\eps,\eta} Du_{\eps,\eta}\in L^2(0,T;L^2),\sqrt\eta\Delta u_{\eps,\eta}\in L^2(0,T;L^2),\\
&\nabla\rho_{\eps,\eta}^{\frac{\gamma}{2}}\in L^2(0,T;L^2),\sqrt\delta \rho_{\eps,\eta}\in L^\infty(0,T;H^{2s+1}),\\
&B_{\eps,\eta}\in L^\infty(0,T;L^2) \cap L^2(0,T;H^1),\sqrt{\delta}\rho_{\eps,\eta}\in L^2(0,T;H^{2s+2}),\\
&\sqrt\kappa\sqrt\rho_{\eps,\eta}\nabla^2\log\rho_{\eps,\eta}\in L^2(0,T;L^2),\sqrt{\nu_b(\rho_{\eps,\eta})}\nabla B_{\eps,\eta}\in L^2(0,T;L^2).
 \end{aligned}
 \end{equation}
 and also satisfies the following
 \begin{equation}\label{e5.17}
\sqrt\kappa\left\|\sqrt\rho_{\eps,\eta}\right\|_{L^2(0,T;H^2)}+\kappa^{1/2}\left\|\nabla\rho_{\eps,\eta}^{1/4}\right\|_{L^4(0,T;L^4)}\leq C.
\end{equation}
where $C$ is independent of the parameters $\eps,\eta,\delta.$

By the above estimate we can get $\nabla\rho_{\eps,\eta}^{-1/2}\in L^2(0,T;L^2)$. 
Since $$\nabla\rho_{\eps,\eta}^{-1/2}=-\frac12\frac{\nabla\rho_{\eps,\eta}}{\rho_{\eps,\eta}^{3/2}}=-\frac{\nabla\sqrt\rho_{\eps,\eta}}{\rho_{\eps,\eta}}.$$
For the case $\rho>1$ then the proof is obvious. then we focus on the case $\rho<1$
\begin{align*}
\nabla\rho^{-1/2}&=\nabla(\rho^{\frac{\gamma^{-}-1}{2}}\rho^{-\frac{\gamma^{-}}{2}})\\
&=\rho^{\frac{\gamma^{-}-1}{2}}\nabla\rho^{-\frac{\gamma^{-}}{2}}+\nabla\rho^{\frac{\gamma^{-}-1}{2}}\rho^{-\frac{\gamma^{-}}{2}}\\
&=\rho^{\frac{\gamma^{-}-1}{2}}\nabla\rho^{-\frac{\gamma^{-}}{2}}+(1-\gamma^{-})\nabla\rho^{-1/2}.
\end{align*}
so $\gamma^{-}\nabla\rho^{-1/2}=\rho^{\frac{\gamma^{-}-1}{2}}\nabla\rho^{-\frac{\gamma^{-}}{2}}$, due to $\int_0^{T}\int_{\O} H_c^{''}(\rho)\left|\nabla\rho\right|^2 dxdt\leq C$ together with the relationship $\rho H_c^{''}(\rho)=P_c^{'}(\rho)$ yield $\int_0^{T}\int_{\O}\rho^{\gamma^{-}-1}\left|\nabla\rho^{-\frac{\gamma^{-}}{2}}\right|^2 dxdt \leq C$.Thus, we have $\nabla\rho^{-1/2}\in L^2(0,T;L^2)$.

Following the same procedure as in the proof of the Lemma 4.1. Application to Aubin-Lions lemma and \eqref{e5.16} and \eqref{e5.17} give rise to the following compactness
\begin{lem}
The following convergence holds
\begin{equation}
\begin{aligned}
&\sqrt\rho_{\eps,\eta}\to \sqrt\rho \,\text{ strongly in} \,L^2(0,T;H^1),\\
&\rho_{\eps,\eta}\to\rho\,\text{strongly in } \,L^2(0,T;H^{2s+1}), \text{and weakly in } L^2(0,T;H^{2s+2}),\\
&\rho_{\eps,\eta} u_{\eps,\eta}\to\rho u \,\text{strongly in}\, L^2(0,T;L^2),\\
&u_{\eps,\eta}\to u\,\text{strongly in}\, L^2(0,T;L^2), B_{\eps,\eta}\to B\, \text{strongly in}\, L^2(0,T;L^2),\\
&\rho_{\eps,\eta}^{-1/2}\to \rho^{-1/2} almost \,everywhere.\nabla B_{\eps,\eta}\to \nabla B\, \text{weakly in}\, L^2(0,T;L^2),\\
&P(\rho_{\eps,\eta})\to P(\rho) \,\text{strongly in }\, L^1((0,T)\times\O),\\
&P_c(\rho_{\eps,\eta})\to P_c(\rho) \,\text{strongly in }\, L^1((0,T)\times\O).
\end{aligned}
\end{equation}
\end{lem}
\begin{pf}
Here we just show the convergence of velocity other convergence can be referred to the previous process. It is noted that the strong convergence of velocity $u$ can be obtained in this step. In fact, since the lower bound of the density is just depend on the parameter $\delta$ and the regularity of $\sqrt\rho\nabla u\in L^2(0,T;L^2)$. Then we get the uniform bound  of $\nabla u\in L^2(0,T;L^2)$ which is independent of $\eps,\eta$. By the strong convergence of $\rho_{\eps,\eta}u_{\eps,\eta}\to\rho u$ we have $\rho_{\eps,\eta}u_{\eps,\eta}\to\rho u$ almost everywhere in $(x,t)\in (0,T)\times \Omega$. together with the convergence of $\rho_{\eps,\eta}^{-1}$ yields that the strong convergence of velocity $u_{\eps,\eta}$. 
\end{pf}
Thus, we can pass to the limit for the nonlinear term $\rho_{\eps,\eta}u_{\eps,\eta}\otimes u_{\eps,\eta}$, Similarly we can also pass to the limit for the capillarity term $\delta\rho_{\eps,\eta}\nabla\Delta\rho_{\eps,\eta}^{2s+1}$, the quantum term and viscosity terms. And we have 
\begin{align*}
&\left|\eta\iint \Delta^2 u_{\eps,\eta}\varphi dx dt\right|\leq \sqrt\eta\left\|\sqrt\eta\Delta u_{\eps,\eta}\right\|_2\left\|\Delta\varphi\right\|_{L^2(0,T;L^2)} \to 0,\\
&\left|\eps\iint\Delta \rho_{\eps,\eta}\varphi dxdt\right|\leq \sqrt\eps\left\|\nabla\rho\right\|_{L^2(0,T;L^2)}\left\|\nabla\varphi\right\|_{L^2(0,T;L^2)}\to 0, \\
&\eps\left|\iint\nabla\rho_{\eps,\eta}\cdot\nabla u_{\eps,\eta} \varphi dxdt\right|\leq C \sqrt\eps\left\|\sqrt\eps\nabla\rho_{\eps,\eta}\right\|_{L^2(0,T;L^2)}\left\|\nabla u_{\eps,\eta}\right\|_{L^2(0,T;L^2)}\to 0.
\end{align*}
Thus we have shown that the limit function $(\rho,u,B)$ is the weak solution of the following system
\begin{equation}
\left\{
\begin{aligned}
&\d \rho+\diver(\rho u)=0, x\in \O, t>0,\\
&\d(\rho u)+\diver(\rho u\otimes u)+\nabla(P(\rho)+P_c(\rho))-2\diver(\rho D(u))\\
&\quad-\delta\rho\nabla\Delta^{2s+1}\rho -2\kappa^2\rho\nabla\left(\frac{\Delta\sqrt\rho}{\sqrt\rho}\right)-(\nabla\times B)\times B=0,\\
&\d B-\nabla\times(u\times B)+\nabla(\nu_b(\rho)\nabla\times B)=0.\\
\end{aligned}
\right.
\end{equation} 
and satisfies the BD entropy estimate \eqref{e5.2} and energy estimate  \eqref{e3.14}. Moreover, we also have the 
\begin{equation}\label{e5.20}
 \kappa\left\|\sqrt\rho\right\|_{L^2(0,T;H^2)}+\kappa^{1/2}\left\|\nabla\rho^{1/4}\right\|_{L^4(0,T;L^4)}\leq C.
\end{equation}
where $C$ is independent of $\delta$.
\subsection{pass to the limit as $\delta\to 0$}
Our goal in this step is to perform the limit as $\delta\to 0$. Here we will lose the uniform lower bound of the density. we can also get additional regularity information. Furthermore the following estimate also holds
\begin{equation}\label{e5.21}
 \begin{aligned}
 &\rho_{\delta}\in L^\infty(0,T;L^\gamma),\rho_{\delta}^{-1}\in L^\infty(0,T;L^{\gamma^{-}}),\nabla\sqrt\rho_{\delta}\in L^\infty(0,T;L^2),\\
&\sqrt\rho_{\delta} u_{\delta}\in L^\infty(0,T;L^2),\sqrt\rho_{\delta} Du_{\delta}\in L^2(0,T;L^2),\\
&\nabla\rho_{\delta}^{\frac{\gamma}{2}}\in L^2(0,T;L^2),\sqrt\delta \rho_{\delta}\in L^\infty(0,T;H^{2s+1}),\\
&B_{\delta}\in L^\infty(0,T;L^2) \cap L^2(0,T;H^1),\sqrt{\delta}\rho_{\delta}\in L^2(0,T;H^{2s+2}),\\
&\sqrt\kappa\sqrt\rho_{\delta}\nabla^2\log\rho_{\delta}\in L^2(0,T;L^2),\sqrt{\nu_b(\rho_{\delta})}\nabla B_{\delta}\in L^2(0,T;L^2).
 \end{aligned}
 \end{equation}
and satisfies the estimate 
\begin{equation}\label{e5.22}
 \kappa\left\|\sqrt\rho_{\delta}\right\|_{L^2(0,T;H^2)}+\kappa^{1/2}\left\|\nabla\rho_{\delta}^{1/4}\right\|_{L^4(0,T;L^4)}\leq C.
\end{equation}
\begin{lem}
From the above estimate \eqref{e5.21} and \eqref{e5.22} , there exists a constant C independent of $\delta$ such that 
\begin{align*}
&\left\|\nabla u\right\|_{L^p(0,T;L^q)}\leq C,\, p=\frac{2\gamma^{-}}{\gamma^{-}+1}, q=\frac{6\gamma^{-}}{3\gamma^{-}+1}, \\
&\left\|u\right\|_{L^P(0,T;L^{q^{\ast}})}\leq C,\, q^{\ast}=\frac{3q}{3-q},\\
&\left\|\sqrt\rho u\right\|_{L^{p^{'}}(0,T;L^{q^{'}})}\leq C, \,p^{'}>2,q^{'}>2.
\end{align*}
\end{lem}
\begin{pf}
 $\nabla u=\frac{1}{\sqrt\rho}\sqrt\rho u$, since $\sqrt\rho u\in L^\infty(0,T;L^2)$ and using the previous estimate we get  $\frac{1}{\sqrt\rho}\in L^{2\gamma^{-}}(0,T;L^{6\gamma^{-}})$, apply the H$\ddot{o}$lder inequality yields that desire estimate. Due to $W^{1,q}$ embedding $L^{q^{\ast}}$, we can obtained $\left\|u\right\|_{L^P(0,T;L^{q^{\ast}})}\leq C$. Next we turn to estimate the $\sqrt\rho u$. First for $0<r<1/2$ we have $\sqrt\rho u=(\sqrt\rho u)^{2r}u^{1-2r}\rho^{1/2-r}$, by the estimate \eqref{e5.21} yield $\sqrt\rho\in L^\infty(0,T;L^3)$ and we get $\rho^{1/2-r}\in L^\infty(0,T;L^{3/(1/2-r)}$. As a consequence we have $(\sqrt\rho u)^{2r}\in L^\infty(0,T;L^{1/r})$ and $u^{1-2r}\in L^{\frac{p}{1-2r}}(0,T;L^{\frac{q^{\ast}}{1-2r}})$
 By  H$\ddot{o}$lder inequality we deduce $\left\|\sqrt\rho u\right\|_{L^{p^{'}}(0,T;L^{q^{'}})}\leq C, \,p^{'}>2,q^{'}>2$ where $\frac{1}{p^{'}}=\frac{1-2r}{p},\frac{1}{q^{'}}=\frac{r}{1}+\frac{1/2-r}{3}+\frac{1-2r}{q^{\ast}}$ then take $\frac{1}{10}<r<\frac12$ can enable us to get the estimate.
\end{pf}

Repeated the same procedure as shown in the above part. we can deduce the following compactness information
\begin{equation}
\begin{aligned}
&\sqrt\rho_{\delta}\to \sqrt\rho \,\text{ strongly in} \,L^2(0,T;H^1),\\
&\sqrt\rho_{\delta} u_{\delta}\to \sqrt\rho u\, \text{strongly in }\, L^2(0,T;L^2),\\
&\rho_{\delta} u_{\delta}\to\rho u \,\text{strongly in}\, L^2(0,T;L^2),\\
&u_{\delta}\to u\,\text{weakly in}\, L^2(0,T;L^2), B_{\delta}\to B\, \text{strongly in}\, L^2(0,T;L^2),\\
&\rho_{\delta}^{-1/2}\to \rho^{-1/2} almost \,everywhere.\nabla B_{\delta}\to \nabla B\, \text{weakly in}\, L^2(0,T;L^2),\\
&P(\rho_{\delta})\to P(\rho) \,\text{strongly in }\, L^1((0,T)\times\O),\\
&P_c(\rho_{\delta})\to P_c(\rho) \,\text{strongly in }\, L^1((0,T)\times\O).
\end{aligned}
\end{equation}

In this step we can get the strong convergence of $\sqrt\rho_{\delta} u_{\delta}$ which is different with the above process. Since from the momentum equation we can get the information of $\d(\rho_{\delta} u_{\delta})$. Apply to the Aubin-Lions lemma yields the almost everywhere convergence of $\rho_{\delta} u_{\delta}$, which, due to the almost everywhere convergence of $\sqrt\rho$. we can get $\sqrt\rho_\delta u_\delta\to \sqrt\rho u$ almost everywhere.Thus we can get the strong convergence of $\sqrt\rho_\delta u_\delta\to \sqrt\rho u$ in $L^2(0,T;L^2)$. The main obstacle in this step is deal with the term $\delta\rho_\delta\nabla\Delta^9\rho_\delta$. Other terms can be easily treated as the above process. 
\begin{lem}
For any test function $\varphi$ we have 
\begin{equation*}
\delta\iint \rho_\delta\nabla\Delta^9\rho_\delta\varphi  dxdt\to 0.
\end{equation*}
\end{lem}
\begin{pf}
By \eqref{e5.21} we have the following uniform estimate
 $$\rho_\delta\in L^\infty(0,T;L^3), \sqrt\delta\rho_\delta\in L^\infty(0,T;H^{2s+1}),\sqrt\delta\rho_\delta\in L^2(0,T;H^{2s+2}).$$ 
Together with the Gagliardo-Nirenberg inequality yields 
\begin{equation}\label{e5.24}
\left\|\nabla^{2s+1}\rho_\delta\right\|_{L^3}\leq C\left\|\rho_\delta\right\|_{W^{2s+2,2}}^{\alpha}\left\|\rho_\delta\right\|_{L^3}^{1-\alpha}.
\end{equation}
where $\alpha\in (0,1)$ satisfy the following relationship 
$$\frac13-\frac{2s+1}{3}=\alpha\left(\frac12-\frac{2s+2}{3}\right)+(1-\alpha)\frac13,\, \alpha=\frac{4s+2}{4s+3}.$$
Moreover, we can also get $\delta^{\frac{\alpha}{2}}\left|\nabla^{2s+1}\rho_\delta\right|\in L^{\frac{\alpha}{2}}(0,T;L^3)$ by using \eqref{e5.24}. Due to
$$\delta\iint \rho_\delta\nabla\Delta^9\rho_\delta\varphi  dxdt= -\delta\iint\Delta^{2s+1}\rho_\delta\Delta^s \diver(\rho_\delta\varphi)dxdt.$$
Then we focus on the most difficult term 
\begin{align*}
&\left|\delta\iint \Delta^{s+1}\rho_\delta\Delta\nabla\rho_\delta\varphi dxdt\right|\\
&\leq C(\varphi) \delta^{\frac{1-\alpha}{2}}\left|\delta^{\frac{\alpha}{2}}\nabla^{2s+1}\rho_\delta\right|_{L^{\frac{2}{\alpha}}(0,T;L^3)}\left|\sqrt \delta \nabla^{2s+2} \rho_\delta\right|_{ L^2(0,T;L^2)}\\
& \to 0 \quad \text{as}\, \delta \to 0.
\end{align*}
Other terms also converge to 0 by the same approach, Thus the proof of the Lemma 6.4 is completed.
\end{pf}
\section{Lower planck limit}
Our goal in this section is to prove the Theorem 3.2. For the sequence of solutions $(\rho_\kappa  ,u_\kappa,B_\kappa)$, then we need  to prove that we can pass to the limit for each term that occurs in the equation. Following the same procedure as shown in section 4 we can get the same compactness information on the solutions. It's worth noting that the strong convergence of $\sqrt\rho_\kappa$ in this step is only in the space $L^2(0,T;L^2)$ rather than $L^2(0,T;H^1)$. Since the uniform bound of $\sqrt\rho_\kappa$ in the space $ L^2(0,T;H^1)$ is just only enable us to get the convergence in $L^2(0,T;L^2)$.
\begin{lem}
Under the condition of Theorem 3.2 and using the uniform estimate we have 
\begin{equation}\label{e6.1}
\begin{aligned}
&\sqrt\rho_{\kappa}\to \sqrt\rho \,\text{ strongly in} \,L^2(0,T;L^2),\\
&\sqrt\rho_{\kappa} u_{\kappa}\to \sqrt\rho u\, \text{strongly in }\, L^2(0,T;L^2),\\
&\rho_{\kappa} u_{\kappa}\to\rho u \,almost\, everywhere,\\
&u_{\kappa}\to u\,\text{weakly in}\, L^2(0,T;L^2),\\
&\rho_{\kappa}^{-1/2}\to \rho^{-1/2} almost \,everywhere,\\
&P(\rho_{\kappa})\to P(\rho) \,\text{strongly in }\, L^1((0,T)\times\O),\\
&P_c(\rho_{\kappa})\to P_c(\rho) \,\text{strongly in }\, L^1((0,T)\times\O),\\
& B_{\kappa}\to B\, \text{strongly in}\, L^2(0,T;L^2),\\
&\sqrt{\nu_b(\rho_\kappa)}\nabla B_{\kappa}\to \sqrt\nu_b(\rho) \nabla B\, \text{weakly in}\, L^2(0,T;L^2).\\
\end{aligned}
\end{equation}
\end{lem}
\begin{pf}
We can use the same procedure as shown in section 3 to obtained the convergence of the above terms. The only different here is the strong convergence of $\sqrt\rho_\kappa$ in the space $L^2(0,T;L^2)$. Due to $\nabla\sqrt\rho\in L^2(0,T;L^2)$ and $\d\sqrt\rho\in L^2(0,T;H^{-1})$, Together with Aubin-Lions lemma yields $\sqrt\rho_\kappa\to \sqrt\rho\, in\, L^2(0,T;L^2)$.
\end{pf}
Next we focus on the quantum term $\rho_\kappa\nabla\left(\frac{\Delta\rho_\kappa}{\sqrt\rho_\kappa}\right)$  
\begin{lem}
For any test function we have 
\begin{equation}
2\kappa^2\iint\rho_\kappa\nabla\left(\frac{\Delta\rho_\kappa}{\sqrt\rho_\kappa}\right) \varphi dxdt\to 0.
\end{equation}
\end{lem}
\begin{pf}
By the inequality \eqref{e4.3} we get
\begin{equation}
\begin{aligned}
&2\kappa^2\left|\iint \rho_\kappa\nabla\left(\frac{\Delta\sqrt\rho_\kappa}{\sqrt\rho_\kappa}\right) \varphi dx dt\right|\\
&\leq 2\kappa^2\left|\iint \sqrt\rho_\kappa\nabla\sqrt\rho_\kappa\nabla\diver\varphi dxdt\right|+4\kappa^2\left|\iint \nabla\sqrt\rho_\kappa\otimes\nabla\sqrt\rho_\kappa\nabla\varphi dxdt\right|.\\
&\leq 2\kappa^2\left\|\sqrt\rho_\kappa\right\|_{L^\infty(0,T;L^2)}\left\|\nabla\sqrt\rho_\kappa\right\|_{L^(0,T;L^2)}\left\|\nabla\diver\varphi\right\|_{L^\infty((0,T)\times \O)}\\
&\quad+4\kappa^2\left\|\nabla\sqrt\rho_\kappa\right\|_{L^(0,T;L^2)}^2\left\|\nabla\varphi\right\|_{L^\infty((0,T)\times \O)}\\
&\leq C\kappa^2 \to 0, \quad as \kappa\to 0.
\end{aligned}
\end{equation}
\end{pf}
Finally, the proof of Theorem 3.2 is completed.

\bibliographystyle{plain}
\bibliography{qmhdreference}
\end{document}